\documentclass[a4paper,12pt, reqno]{amsart}
\usepackage{a4,wasysym}
\usepackage{amsmath,graphicx, amssymb,color,ulem,bbm}

\newtheorem{theorem}{Theorem}
\newtheorem{corollary}[theorem]{Corollary}

\newtheorem{proposition}[theorem]{Proposition}

\newcommand{\ud}{\mathrm{d}}

\newcommand{\N}{\mathcal{N}}

\newcommand{\eofproof}{\hfill$\Box$ }

\begin{document}

\title[On a strain-structured epidemic model]{On a strain-structured epidemic model}

\author[\`{A}. Calsina]{\`{A}ngel Calsina}
\address{\`{A}ngel Calsina, Department of Mathematics, Universitat Aut\`{o}noma de Barcelona, Bellaterra, 08193, Spain}
\email{acalsina@mat.uab.es}

\author[J\'{o}zsef Z. Farkas]{J\'{o}zsef Z. Farkas}
\address{J\'{o}zsef Z. Farkas, Division of Computing Science and Mathematics, University of Stirling, Stirling, FK9 4LA, Scotland, United Kingdom }
\email{jozsef.farkas@stir.ac.uk}

\subjclass{ }
\date{\today}

\begin{abstract}
\begin{sloppypar}
We introduce and investigate an SIS-type model for the spread of an infectious disease, where the infected population is structured with respect to the different strain of the virus/bacteria they are carrying. Our aim is to capture the interesting scenario when individuals infected with different strains cause secondary (new) infections at different rates. Therefore, we consider a nonlinear infection process, which generalises the bilinear process arising from the classic mass-action assumption. Our main motivation is to study competition between different strains of a virus/bacteria. From the mathematical point of view, we are interested whether the nonlinear infection process leads to a well-posed model. We use a semilinear formulation to show global existence and positivity of solutions up to a critical value of the exponent in the nonlinearity. Furthermore, we establish the existence of the endemic steady state for particular classes of nonlinearities.   \end{sloppypar}
\end{abstract}
\maketitle

\section{Introduction}

Mathematical modelling of epidemiological processes, such as an influenza outbreak, frequently involves differential equations. 
There is a vast literature of ordinary differential equations used to model infection dynamics of various diseases, see e.g. \cite{DH}. In recent years the emphasis on modelling has drifted to take into account the complex network structure of susceptible and infected populations. The network structure can naturally account for changes in the infection dynamics, see e.g. the recent papers \cite{Green,Kiss1,Kiss2} and the references therein. In this context, and motivated by recent disease outbreaks, the effect of so called super-spreading individuals on the infection dynamics was analysed. We may say that super-spreaders are infected individuals that cause larger than average number of secondary infections. This phenomenon has been observed in various diseases, for example in case of SARS, see e.g. \cite{LlSm}. In the first instance it may be natural to assume that super-spreaders are individuals who, due to the complex network structure of the population, encounter more contacts with susceptible individuals on average. However, there are different possible causes why some individuals may transmit the infection with higher probability. For example, an individual may be a more suitable host (acting as an environment) for a particular strain of the virus/bacteria. Hence virus/bacterium replication dynamics will depend on the individual (host). As a result, individuals carrying significantly larger quantities of the virus/bacteria may become more infectious. Another, and perhaps more interesting scenario is when simply, due for example to (rapid) evolution, the virus/bacteria exhibit a large number of different genetic strains. Individuals carrying different strains may then cause significantly different number of secondary infections. 

It has been also increasingly recognised that the modelling of infectious diseases may require structured models. 
It is often the case that individuals may need to be distinguished according to an inherent property. 
Several structured population models were formulated to model the infection dynamics of a disease, see for example \cite{CF,Ducrot,FH}. In earlier models structuring of infected individuals was often based on infection-age, see e.g. \cite{FPW}. More recently in \cite{CF,FH} structured models were considered with Wentzell (or Feller) boundary conditions. These boundary conditions provide an elegant way to formulate a structured SIS-type model, where the uninfected/susceptible individuals correspond to the boundary state. A very similar idea motivated the early paper \cite{WHG}.  Still, structured population models (formulated as partial differential equations) are employed less frequently than their unstructured counterparts. This is mainly because of the challenges arising in the analysis of those models. 
 
Here we consider a fairly straightforward SIS-type model. We model a population which is infected with a virus/bacteria which has a large number of different strains. Hence we choose to model the type of strain by a continuous variable, see e.g. \cite{KD}. We assume that different strains may have different infective capabilities. That is, an infected individual with strain $y$ will cause a number of secondary infections depending on the value $\gamma(y)$ of a function $\gamma$. Typically, one may assume that $\gamma$ is a monotone increasing function, but from the mathematical point of view we do not have to impose such restriction. The classic mass-action assumption, which usually leads to a bilinear model is modified to take into account the highly nonlinear effect of super-spreaders on the dynamics. This, combined with structuring of the infected population, leads to an infinite dimensional nonlinearity in the (non-local) recruitment operator.  
Different nonlinear incidence rates were introduced and studied in the context of ordinary differential equation models 
already in the seminal papers \cite{CAP1,CS}, and more recently for example in \cite{KOR}. In \cite{CS} the bilinear incidence rate was modified to take into account saturation effects for new infections for large population sizes of the infected individuals. In the more recent paper \cite{KOR} ordinary differential equations with very general nonlinear transmission terms were considered. The emphasis was on constructing Lyapunov functions for the models, and on numerical simulations.

Our model in principle allows to study competition between different strains of the virus/bacteria. For example, the question whether a particular strain could win the evolutionary race can be answered by establishing global existence of solutions (or the lack of it) of the partial differential equation. Blow-up at a single value of the structuring variable may be interpreted that a single strain of the bacteria/virus wins the evolutionary race. First we will establish global existence of solutions when $\gamma$ is below a critical value. Then we will study the existence  of positive steady states for different values of $\gamma$. In particular, we will establish the existence of endemic steady states for the bilinear and for the quadratic 
case. Our results closely resemble the corresponding results for the counterpart ordinary differential equation model.

\section{The model}

We consider the fundamental scenario when individuals in the population are grouped into two classes: susceptible and infected, with the infected individuals structured according to the particular strain of the virus/bacteria they carry. The type of strain is denoted by the variable $x$. For simplicity we assume that $x$ takes values from the interval $[0,1]$.  We denote by $v$ the density of infected individuals, and by $S$ the number of susceptible individuals. 
We consider the following model.
\begin{equation}\label{superspread}
\begin{aligned}
v_t(x,t)-\left(d(x)v_x(x,t)\right)_x &= -\varrho(x)v(x,t)+S(t)\int_0^1\beta(x,y)v(y,t)^{1+\gamma(y)}\,\ud y, \\
S'(t)=\int_0^1\varrho(x)v(x,t)\,\ud x&-S(t)\int_0^1\int_0^1\beta(x,y)v(y,t)^{1+\gamma(y)}\,\ud y\,\ud x, \\
v_x(0,t)&=v_x(1,t)=0,\quad v(x,0)=v_0(x),\quad S(0)=S_0.
\end{aligned}
\end{equation}
In the model above $\varrho$ denotes the recovery rate, that is the rate of return of infected indi\-vi\-duals into the susceptible class. 
The integral term in the partial differential equation above describes the infection process. An infected individual carrying (mainly) strain $y$, 
transmits the bacteria/virus upon contact to a susceptible individual. Due to the fact that the newly infected individual represents a (possibly)  different environment, a possibly (different) strain $x$ is (instantaneously) selected. From the mathematical point of view this is captured by the function $\beta$, which is the rate at which individuals of strain $y$ produce newly infected individuals of strain $x$.  
As we mentioned in the introduction, we assume that individuals carrying different strains (may) have different infectiousness. 
This is captured by the function $\gamma$.  Apart from the infection process, changes of the strain/infectiousness due to random mutations inside the host are modelled by the diffusion rate $d$.

In the rest of the paper we will also assume that the model ingredients satisfy the following natural assumptions.
\begin{align}
& 0\le\beta\in C^1([0,1]^2),\,\,||\beta||_{C^1}\le b,\quad 0<d_0\le d\in C^1(0,1),\,\, ||d||_{C^1}\le d_1, \nonumber \\
& 0\le\varrho\in C^1(0,1),\,\, ||\varrho||_{C^1}\le r,\quad 0\le\gamma\le\Gamma. \label{hypo}
\end{align}

Note that integrating the first equation of \eqref{superspread} and adding it to the second one,  
and applying the boundary conditions, one obtains 
\begin{equation*}
\frac{\ud}{\ud t} \left(\int_0^1v(x,t)\,\ud x+S(t)\right)=0.
\end{equation*}
Hence the total population size of any solution with initial condition $(v_0,S_0)$ is preserved, i.e. it remains $\int_0^1 v_0(x)\,\ud x+S_0$ for all times. Here we concentrate on the infection dy\-na\-mics, which for a number of diseases may be assumed to take place on a faster time-scale.  That is we neglect population dynamics. We note that it would be straightforward to modify model \eqref{superspread} to incorporate birth and death processes; and also to incorporate additional compartments  for recovered (R) and for exposed (E) individuals. 

One may naturally think of model \eqref{superspread} as a continuum approximation of the following system of ODEs, for an infection 
with a finite number of strains. 
\begin{align}
& \frac{\ud}{\ud t}I_i(t)=\sum_{j=1}^n\left(d_{i,j}I_j(t)-d_{j,i}I_i(t)\right)-\varrho_i I_i(t)+S(t)\sum_{j=1}^n\beta_{i,j}I_j^{1+\gamma_j}(t), \quad i=1,\cdots,n,\nonumber \\
& \frac{\ud}{\ud t}S(t)=\sum_{i=1}^n\varrho_i I_i(t)-S(t)\sum_{i=1}^{n}\sum_{j=1}^{n}\beta_{i,j} I_j^{1+\gamma_j}(t). \label{ODEmodel-0}
\end{align}
In the simplest case of only one infection strain, model \eqref{ODEmodel-0} reduces to
\begin{align}
& \frac{\ud}{\ud t}I(t)=-\varrho I(t)+\beta S(t)I^{1+\gamma}(t), \nonumber \\
& \frac{\ud}{\ud t}S(t)=\varrho I(t)-\beta S(t)I^{1+\gamma}(t). \label{ODEmodel}
\end{align}
Model \eqref{ODEmodel} is globally well posed for positive initial conditions and for any $\gamma\ge 0$, as the total population size is preserved. 
This may not be the case for the structured model \eqref{superspread}. In particular, we conjecture that for $\gamma>1$ blow-up might be possible on a set of measure zero. We will return later again to compare the qualitative properties of the partial differential equation model \eqref{superspread} to the unstructured (counterpart) model \eqref{ODEmodel}.

\section{Existence and positivity of solutions}

In this section we investigate the existence and positivity of solutions of model \eqref{superspread}. 
Our goal is to use a semilinear formulation based on the variation of constants formula, as in \cite{Henry}. 
One of the advantages of this formulation is that the principle of linearised stability can be established.  
We note that the question of (both) local and global existence are far from trivial due to the superlinear infection process. 
There is a vast literature of non-existence results for nonlinear reaction-diffusion equations. 
Here we only mention the paper \cite{SOU} for further reference.   

\subsection{Local existence and positivity}
First we establish local existence and positivity of solutions of \eqref{superspread} for $\gamma\ge 0$. To this end, we recast model \eqref{superspread} in the form of a Cauchy problem as
\begin{equation}\label{ACauchy}
\frac{\ud {\bf v}}{\ud t}=\mathcal{L}{\bf v}+\mathcal{N}{\bf v}, \quad t>0,\quad\quad {\bf v}(0)={\bf v}_0.  
\end{equation}
Above we introduced the vector notation ${\bf v}=(v,S)$. The linear operator $\mathcal{L}$ on the state space  $\mathcal{X}=L^1(0,1)\times\mathbb{R}$, is naturally defined as the following closed unbounded linear operator
\begin{equation}\label{linop}
\begin{aligned}
\mathcal{L}{\bf v}&=
\begin{pmatrix}
\frac{\partial}{\partial x}\left(d(\cdot)\frac{\partial v}{\partial x}\right)-\varrho(\cdot)v \\
 \int_0^1\varrho(x)v(x)\,\ud x\\
\end{pmatrix},\\
\text{D}(\mathcal{L})&=  \left\{ (v,S)^T\in W^{2,1}(0,1)\times\mathbb{R} \: :\: v'(0)=v'(1)=0\right\}.
\end{aligned}
\end{equation}
Note that $\mathcal{L}$ generates an analytic semigroup on the biologically relevant state space $\mathcal{X}=L^1(0,1)\times\mathbb{R}$. The nonlinear operator $\mathcal{N}$ is defined as
\begin{equation}\label{nonlinop}
\mathcal{N}{\bf v}=\begin{pmatrix} S\int_0^1\beta(\,\cdot\,,y)|v(y)|^{1+\gamma(y)}\,\ud y\\
-S\int_0^1\int_0^1\beta(x,y)|v(y)|^{1+\gamma(y)}\,\ud y\,\ud x\\
\end{pmatrix}.
\end{equation}
Above in the definition of $\mathcal{N}$ we have taken the absolute value of $v$ in order for $\mathcal{N}$ being well-defined on linear spaces, that is also for functions which are not in the positive cone. This is not a real constraint, as we are only interested in positive solutions of \eqref{superspread}. The nonlinear operator $\mathcal{N}$ cannot be defined on the whole space $\mathcal{X}$ because 
$\gamma\ge 0$, (unless $\gamma\equiv 0$, which corresponds to the classic bilinear infection process). Hence to establish existence and uniqueness of solutions of \eqref{superspread} we use the framework of fractional power spaces. The method developed in \cite{Henry} utilises some intermediate spaces between the space $\mathcal{X}$ and the domain of the generator of the linear part $\mathcal{L}$. We also note that the theory developed in \cite{Henry} relies heavily on the fact that the linear operator $\mathcal{L}$ generates an analytic semigroup. 

Let us recall from \cite{Henry} that if $\mathcal{A}$ is a sectorial operator with spectral bound $s(-\mathcal{A})<a$, then for $\mathcal{A}_1:=\mathcal{A}+a\,\mathcal{I}$ and any $0<\alpha$ we can define 
\begin{equation}
\mathcal{A}_1^{-\alpha}=\frac{1}{\Gamma(\alpha)}\int_0^\infty t^{\alpha-1}\exp(-\mathcal{A}_1t)\, \ud t,
\end{equation}
where $\exp\left(-\mathcal{A}_1 t\right)$ is the analytic semigroup generated by $-\mathcal{A}_1$ and $\Gamma$ stands for the gamma function. $\mathcal{A}_1^\alpha$ is then defined as the inverse of $\mathcal{A}_1^{-\alpha}$. 
Moreover, the fractional power space $\mathcal{X}^\alpha$ is then defined for $\alpha>0$ as 
\begin{equation*}
\mathcal{X}^\alpha=\text{D}\left(\mathcal{A}_1^\alpha\right),\quad  \text{with norm}\quad ||x||_{\mathcal{X}^\alpha}=\left|\left|\mathcal{A}_1^\alpha\, x\right|\right|.
\end{equation*} 
In our setting, choosing $\mathcal{A}=-\mathcal{L}$, we have:
\begin{equation*}
\text{D}(\mathcal{L})=\mathcal{X}^1\subset\mathcal{X}^{\alpha}\subset\mathcal{X}^0=L^1(0,1)\times\mathbb{R}=\mathcal{X}, \quad \alpha\in (0,1).
\end{equation*} 
One of the cornerstones of the applicability of the framework developed in \cite{Henry} is the choice of an appropriate and convenient intermediate fractional power space $\mathcal{X}^\alpha$. We start by recalling Theorem 1.6.1 from \cite{Henry} for the reader's convenience.
\begin{theorem}(Theorem 1.6.1, \cite{Henry})\label{embeddingth}
Assume that $\Omega\subset\mathbb{R}^n$ is an open set with smooth boundary. Furthermore assume that $1\le p<\infty$, and $\mathcal{A}$ is a sectorial 
operator in $\mathcal{X}=L^p(\Omega)$ with $D(\mathcal{A})\subset W^{m,p}(\Omega)$, for some $m\ge 1$. Then for $0\le\alpha\le 1$
\begin{align}
& \mathcal{X}^\alpha\subset W^{k,q}(\Omega),\quad \text{if}\quad k-\frac{n}{q}<\alpha m-\frac{n}{p} \quad \text{and}\quad p\le q, \\
& \mathcal{X}^\alpha\subset C^\nu(\Omega),\quad \text{if}\quad 0\le\nu<m\alpha-\frac{n}{p}.
\end{align}
\end{theorem}

In our setting we have $\mathcal{X}=L^1(0,1)\times\mathbb{R}$, and $-\mathcal{L}$ is a sectorial operator with two components, and with domain defined in \eqref{linop}. $\mathcal{X}^\alpha$ is the range of $\left(-\mathcal{L}+a\mathcal{I}\right)^{-\alpha}$. Since the second component of $(-\mathcal{L}+a\mathcal{I})^{-\alpha}$ is scalar, it follows that its range $\mathcal{X}^\alpha=\left(L^1(0,1)\right)^\alpha\times\mathbb{R}$. Theorem \ref{embeddingth} in particular implies that
\begin{equation}\label{interspace}
\mathcal{X}^\alpha=\left(L^1(0,1)\right)^\alpha\times\mathbb{R}\subset W^{1,1}(0,1)\times\mathbb{R},\quad \text{for}\quad \frac{1}{2}<\alpha\le 1.
\end{equation}
At the same time we note that the largest space on which the nonlinearity is readily defined is $W^{1,1}(0,1)\times\mathbb{R}$, since its elements are the absolutely continuous functions. On the other hand, for $\alpha<\frac{1}{2}$ the fractional power spaces $\mathcal{X}^\alpha$ contain unbounded functions.  Hence we will work using the intermediate space $\mathcal{X}^\alpha$,  where $\alpha\in\left(\frac{1}{2},1\right)$. 
Also note that $W^{1,1}(0,1)\subset L^\infty(0,1)$, with continuous inclusion. 

In the rest of the section we will use the notation $||\cdot||_1$ to denote the usual 
norm on $L^1(0,1)$, while $||\cdot||_{W^{1,1}}$ will denote the usual norm on $W^{1,1}(0,1)$. 
\begin{proposition} The operator $\N$ maps $\mathcal{X}^{\alpha}$, for any $\alpha\in\left(\frac{1}{2},1\right)$, into $\mathcal{X}=L^1(0,1)\times\mathbb{R}$, and it is locally Lipschitz continuous. 
\end{proposition}
{\bf Proof.}
We need to show that if ${\bf u}\in\mathcal{X}^\alpha$ then there exists a neighbourhood $U$ of ${\bf u}$, and a real number $C$,   such that for any ${\bf v,w}\in U\Rightarrow||\mathcal{N}({\bf v})-\mathcal{N}({\bf w})||_{\mathcal{X}}\le C\,||{\bf v}-{\bf w}||_{\mathcal{X}^{\alpha}}$. We have the following estimate:
\begin{align}
||\mathcal{N}({\bf v})-\mathcal{N}({\bf w})||_{\mathcal{X}} \le & \,2\left|\left|S_v\int_0^1\beta(\cdot,y)|v(y)|^{1+\gamma(y)}\,\ud y-S_w\int_0^1\beta(\cdot,y)|w(y)|^{1+\gamma(y)}\,\ud y\right|\right|_1 \nonumber \\
 \le & \,2\,|S_v|\,\left|\left|\int_0^1\beta(\cdot,y)\left(|v(y)|^{1+\gamma(y)}-|w(y)|^{1+\gamma(y)}\right)\,\ud y\right|\right|_1 
\label{estimate1} \\
& +2\left|S_v-S_w\right|\left|\left|\int_0^1\beta(\cdot,y)|w(y)|^{1+\gamma(y)}\,\ud y\right|\right|_1. \label{estimate2}
\end{align}
Above, $S_v$ and $S_w$ stand for the second component of the vector ${\bf v}$ and ${\bf w}$, respectively. For any $y\in [0,1]$ we have
\begin{align*}
\left||v(y)|^{1+\gamma(y)}-|w(y)|^{1+\gamma(y)}\right|=\left||v(y)|-|w(y)|\right|(1+\gamma(y))\, z^{\gamma(y)}, 
\end{align*}
where $z\in(|v(y)|,|w(y)|)$.
Note that if for some $y$ we have $|v(y)|=|w(y)|$, then the left-hand side of the equality above equals zero. 
We have
\begin{align*}
\left||v(y)|^{1+\gamma(y)}-|w(y)|^{1+\gamma(y)}\right|=||v(y)|-|w(y)||(1+\gamma(y))\max\left\{|v(y)|^{\gamma(y)},|w(y)|^{\gamma(y)}\right\}.
\end{align*}
With this, we obtain the following upper bound for the expression in \eqref{estimate1}
\begin{align*}
& |S_v|\,||\beta||_\infty(1+\Gamma)\,\max\left\{1+\sup_{y\in [0,1]}|v(y)|^\Gamma,1+\sup_{y\in [0,1]}|w(y)|^\Gamma\right\}\,\int_0^1|v(y)-w(y)|\,\ud y \\
& \le |S_v|||\beta||_\infty(1+\Gamma)\,c\,\max\left\{1+||v||_{W^{1,1}}^\Gamma,1+||w||_{W^{1,1}}^\Gamma\right\}\,||v-w||_1 \\
& \le C_1\, ||v-w||_{W^{1,1}}\le C_1||{\bf v}-{\bf w}||_{W^{1,1}(\Omega)},
\end{align*}
on any bounded set $U\subset W^{1,1}(\Omega)$. Similarly, we obtain the following upper bound for the expression in \eqref{estimate2}
\begin{align*}
& |S_v-S_w|||\beta||_\infty\int_0^1 \max\left\{1,||w||_\infty^{1+\Gamma}\right\}\,\ud y\le |S_v-S_w|||\beta||_\infty\left(1+||w||_\infty^{1+\Gamma}\right) \\
& \le |S_v-S_w|||\beta||_\infty \left(1+c||w||_{W^{1,1}}^{1+\Gamma}\right)\le C_2 |S_v-S_w|\le C_2||{\bf v}-{\bf w}||_{W^{1,1}(\Omega)}.
\end{align*}
Since the inclusion of $\mathcal{X}^\alpha$ for $\alpha\in\left(\frac{1}{2},1\right)$ in Theorem \ref{embeddingth} is continuous, we have showed that $\mathcal{N}\,:\,\mathcal{X}^\alpha\to\mathcal{X}$ is locally Lipschitz continuous for $\alpha\in\left(\frac{1}{2},1\right)$. \eofproof

Theorem 3.3.3 in \cite{Henry} guarantees the existence of a local solution, which we formally state in the following corollary.

\begin{corollary}\label{localex}
For any $(v_0,S_0)^T={\bf v}_0\in\mathcal{X}^\alpha,\,\alpha\in\left(\frac{1}{2},1\right)$ there exists a $t_1>0$ such that model \eqref{superspread} 
has a unique solution ${\bf v}(t)$ for $t\in[0,t_1)$.
\end{corollary}
Next we establish positivity of the local solution. 

\begin{proposition}
For ${\bf v}_0\in\mathcal{X}_+^\alpha,\,\alpha\in\left(\frac{1}{2},1\right)$, the solution ${\bf v}(t)$ of model \eqref{superspread} is positive for the time of its existence.
\end{proposition}
{\bf Proof.}
We rewrite the abstract Cauchy problem \eqref{ACauchy} as follows:
\begin{equation}\label{newCauchy}
\frac{\ud {\bf v}}{\ud t}=\mathcal{L}{\bf v}-\kappa{\bf v}+\mathcal{N}{\bf v}+\kappa{\bf v}, \quad t>0,\quad {\bf v}(0)={\bf v}_0,   
\end{equation}
where $\kappa$ is a positive real number, to be chosen later. $\mathcal{L}$ and $\mathcal{N}$ are defined earlier in \eqref{linop} and \eqref{nonlinop}, respectively. 
Note that $\mathcal{L}-\kappa\,\text{I}$ is a generator of an analytic and positive semigroup $\mathcal{T}_\kappa(t)$, for any $\kappa>0$. 

For any given initial condition ${\bf v}_0\in\mathcal{X}^\alpha_+$ we can choose an open ball ${\bf B}$ of radius $B$ that contains ${\bf v}_0$. On the ball ${\bf B}$ we have
\begin{align}
\int_0^1\int_0^1\beta(x,y)|v(y)|^{1+\gamma(y)}\,\ud y\,\ud x & \le \,||\beta||_\infty
\int_0^1 \max\{1,\sup_{y\in [0,1]}v(y)^{1+\Gamma}\}\,\ud y \nonumber \\
& \le \, b\left(1+c_0||v||_{W^{1,1}(0,1)}^{1+\Gamma}\right)\le \, b\left(1+c_1||{\bf v}||_{\mathcal{X}^\alpha}^{1+\Gamma}\right) \nonumber \\
& \le b\left(1+c_1B^{1+\Gamma}\right)=:\kappa(B)=\kappa.
\end{align}
This shows that the second component of the nonlinear map $\mathcal{N}+\kappa\,\text{I}$ is positive on the ball of radius $B$ intersected with the positive cone of $\mathcal{X}^\alpha$, whereas the first component of $\mathcal{N}$ is positive on the whole positive cone of $\mathcal{X}^\alpha$. 

Also, there exists a set ${\bf B}_\delta$ contained in ${\bf B}$ such that ${\bf B}_\delta\,\cap\,\mathcal{X}^\alpha_+$ is closed. 

As in the proof of Theorem 3.3.3 in \cite{Henry} the solution of \eqref{newCauchy} can be obtained by the contraction mapping principle, that is by the following iteration
\begin{equation}\label{iterate}
{\bf v}_{n+1}(t)=\mathcal{T}_\kappa(t){\bf v}_0+\int_0^t\mathcal{T}_\kappa(t-s)\left(\mathcal{N}+\kappa\,\text{I}\right)({\bf v}_n(s))\,\ud s,\quad {\bf v}_0(t)={\bf v}_0\in\mathcal{X}_+^\alpha,\quad t\in [0,t_*).
\end{equation}
For small enough $t_*$ the right hand side of \eqref{iterate} leaves the set $C\left([0,t_*],{\bf B}_\delta\right)$ invariant, which implies that 
the sequence is positive, and therefore its limit is also positive. 

Finally, let us assume that the solution ${\bf v}(t)$ is not positive 
for every $t\in [0,t_1)$. Let 
\begin{equation*}
T_*=\sup_{t\in[0,t_1)}\{t\, |\, {\bf v}(s)\in\mathcal{X}^\alpha_+\,\text{for all}\, s\in [0,t]\}.
\end{equation*}
Since $0<T_*<t_1$, and ${\bf v}$ is a continuous function, we have 
\begin{equation*}
\displaystyle\lim_{s\to T_*^-}{\bf v}(s)={\bf v}(T_*)\in\mathcal{X}^\alpha_+.
\end{equation*} 
Taking ${\bf v}(T_*)$ as the initial condition in the Cauchy problem \eqref{newCauchy} above, and applying the arguments of the proof in the previous paragraph, we get a contradiction. \eofproof

\subsection{Global existence}

Next we establish the existence of a unique non-negative global solution of \eqref{superspread}, for the case $\Gamma\le 1$. 
\begin{theorem}\label{global}
For $\Gamma\le 1$ and ${\bf v}_0\in\mathcal{X}_+^\alpha,\,\alpha\in\left(\frac{1}{2},1\right)$, model \eqref{superspread} admits a unique non-negative solution for 
all $t\ge 0$.
\end{theorem}
{\bf Proof.} First note that, similarly to \eqref{estimate1}-\eqref{estimate2}, we obtain 
\begin{align}
||\mathcal{N}({\bf v})||_1 & \le 2\,|S|\int_0^1\beta(x,y)\int_0^1 |v(y)|^{1+\gamma(y)}\,\ud y\,\ud x \nonumber \\
 & \le 2\,|S|\, ||\beta||_\infty \int_0^1 \displaystyle\max\{1, |v(y)|^{1+\Gamma}\}\,\ud y \nonumber \\
 & \le 2\,b\,|S| \left(1+\displaystyle\sup_{y\in [0,1]} |v(y)|^{1+\Gamma}\right) \nonumber \\
 & \le 2\,b\,|S| \left(1+c_0||v||_{W^{1,1}}^{1+\Gamma}\right). \label{glest1}
\end{align}
This shows that $\mathcal{N}$ maps bounded subsets of $\mathcal{X}^{\alpha}$ for $\alpha\in\left(\frac{1}{2},1\right)$ into bounded sets of $\mathcal{X}$. Theorem 3.3.4 in \cite{Henry} states that if $\mathcal{N}$ maps bounded sets into bounded sets, then either the maximal time of existence $t_1$, is infinite or $||{\bf v}(t)||_{\mathcal{X}^{\alpha}}\to \infty$ as $t\to t_1$. Therefore, to establish global existence, it suffices to obtain an a priori bound for $||{\bf v}(t)||_{\mathcal{X}^{\alpha}}$ 
for $t\in (0,t_1)$. This will be done in three steps. We will begin by obtaining an a priori bound of the solution in $W^{1,1}(\Omega)$. Then we will use this to obtain an a priori bound of the nonlinearity along the local solution. Finally, we use the variation of constants formula to obtain the needed a priori bound on the $\mathcal{X}^\alpha$ norm of the solution.

Since the $L^1$ norm of the solution ${\bf v}$ is bounded by the (constant) total population size  $P_*=\int_0^1 v_0(x)\,\ud x+S_0$,  obtaining 
the $W^{1,1}(\Omega)$ bound amounts to prove that the derivative $v_x$ of the first component of the solution ${\bf v}$ remains bounded in $L^1$, while $t$ approaches the upper bound of the existence time interval.

The time derivative of the local solution is actually in $\mathcal{X}^{\alpha}$ by Theorem 3.5.2, hence we can differentiate the first equation in \eqref{superspread} to obtain for $u(x,t)=v_x(x,t)$ an equation, which is well-defined in $L^1$.
\begin{align}
 u_t(x,t)-\left(d(x)u(x,t)\right)_{xx}=& -\varrho(x)u(x,t)-\varrho'(x)v(x,t) \nonumber \\
 & +\left(P_*-\int_0^1v(x,t)\,\ud x\right)\int_0^1\beta_x(x,y)v(y,t)^{1+\gamma(y)}\,\ud y, \nonumber \\
u(x,0)&=u_0(x)=v_0'(x),\quad x\in [0,1], \nonumber \\
 u(0,t)&=u(1,t)=0, \quad  t\in(0,t_1).\label{eq-diff}
\end{align}
The first equation of \eqref{eq-diff} is a linear inhomogeneous equation in $u$, and the solution for $t\in (0,t_1)$ (which is the time of existence of the local solution $v$) can be written as:
\begin{equation}
u(t)=\mathcal{U}(t)u_0+\int_0^t\mathcal{U}(t-s)f(s)\,\ud s,\quad t\in (0,t_1),
\end{equation} 
where 
\begin{equation}
f(t)=-\varrho'(\cdot)v(\cdot,t)+\left(P_*-\int_0^1v(x,t)\,\ud x\right)\int_0^1\beta_x(\cdot,y)v(y,t)^{1+\gamma(y)}\,\ud y,\quad  t\in (0,t_1),
\end{equation}
and $\mathcal{U}(t)$ is the analytic and positive semigroup generated by the linear operator $\mathcal{B}$, 
which is defined as
\begin{align*}
& \mathcal{B}\, u= \left(d(\cdot)u(\cdot)\right)_{xx}-\varrho(\cdot)u,\quad \text{D}(\mathcal{B})=\{u\in W^{2,1}(0,1)\,|\, u(0)=u(1)=0\}.
\end{align*}
We have
\begin{align}
||u(t)||_1 & \le ||\mathcal{U}(t)u_0||_1+\left|\left|\int_0^t\mathcal{U}(t-s)f(s)\,\ud s\right|\right|_1 \nonumber \\
& \le C_1\exp(\omega t)||u_0||_1+\int_0^t C_1\exp(\omega(t-s))\,||f(s)||_1\,\ud s. \label{u1estimate}
\end{align}
The assumptions on $\beta$ and $\varrho$ in \eqref{hypo} allow us to obtain the following estimate.
\begin{equation*}
||f(s)||_1\le P_*\left(r+b\left(1+||v(s)||_{1+\Gamma}^{1+\Gamma}\right)\right),\quad s\in(0,t_1).
\end{equation*}
Hence we need to obtain an a priori $L^{1+\Gamma}$ bound for the local solution $v$. Note that the local solution $v$ belongs to the domain of $\mathcal{L}$ and satisfies
\begin{equation}\label{equationv}
v_t(x,t)-(d(x)v_x(x,t))_x=-\varrho v(x,t)+\left(P_*-\int_0^1v(y,t)\,\ud y\right)\int_0^1\beta(x,y)v(y,t)^{1+\gamma(y)}\,\ud y.
\end{equation} 
We multiply equation \eqref{equationv} by $v^\Gamma$ and integrate from $0$ to $1$ to obtain
\begin{align}
\frac{\ud}{\ud t}\left(\frac{\int_0^1 v^{1+\Gamma}(x,t)\,\ud x}{1+\Gamma}\right)& =-\int_0^1 d(x)v^2_x(x,t)\Gamma v^{\Gamma-1}(x,t)\,\ud x-\int_0^1\varrho(x)v^{1+\Gamma}(x,t)\,\ud x \nonumber \\
& +\left(P_*-\int_0^1 v(y,t)\,\ud y\right)\int_0^1 v^\Gamma(x,t)\int_0^1\beta(x,y)v(y,t)^{1+\gamma(y)}\,\ud y\,\ud x. 
\nonumber \\
& \Rightarrow  \\
\frac{\ud }{\ud t}\left(||v(t)||_{1+\Gamma}^{1+\Gamma}\right) & \le (1+\Gamma)P_* \,b\, ||v(t)||_\Gamma^\Gamma\int_0^1 \left|v(y,t)^{1+\gamma(y)}\right|\,\ud y  \nonumber \\
& \le (1+\Gamma)P_*\, b\,  \left(1+||v(t)||_1\right) \left(1+||v(t)||_{1+\Gamma}^{1+\Gamma}\right) \nonumber \\
& \le (1+\Gamma)P_*\,b\,(1+P_*)\left(1+||v(t)||_{1+\Gamma}^{1+\Gamma}\right),\quad t\in (0,t_1), \label{intparts}
\end{align}
for $\Gamma\le 1$ and $v\ge 0$. Note that if the first component of the solution ${\bf v}$ vanishes at some point(s), the estimates above can 
still be obtained by integrating (instead on $(0,1)$) separately on the union of the disjoint intervals where $v$ does not vanish. 

So for $t\in (0,t_1)$, which is the maximal time interval of local existence, we obtained that $||v(t)||_{1+\Gamma}^{1+\Gamma}$ satisfies the ordinary differential inequality \eqref{intparts}. 
That is, there exists a constant $k_1$ such that 
\begin{equation}
||v(t)||_{1+\Gamma}^{1+\Gamma}\le k_1\,\exp((1+\Gamma)P_*\,b\,(1+P_*)\,t).
\end{equation}
Hence there exists a constant $k_2$ such that
\begin{equation}
||f(t)||_1\le k_2\exp((1+\Gamma)P_*\,b\,(1+P_*)\,t),\quad t\in (0,t_1).
\end{equation}
From equation \eqref{u1estimate} we obtain:
\begin{equation}
||u(t)||_1\le  C_1\exp(\omega t)||v_0||_{W^{1,1}(0,1)}+C_2\exp\left\{(1+\Gamma)P_*\,b\,(1+P_*)\,t\right\},\quad t\in (0,t_1).
\end{equation}
This implies the needed a priori bound on $||v(t)||_{W^{1,1}(0,1)}$ as follows.
\begin{equation}
||v(t)||_{W^{1,1}(0,1)}\le P_*+ C_1\exp(\omega t)||v_0||_{W^{1,1}(0,1)}+C_2\exp\left\{(1+\Gamma)P_*\,b\,(1+P_*)\,t\right\},\quad t\in (0,t_1).
\end{equation}
This, together with \eqref{glest1}, implies that for any initial condition there exists a conti\-nu\-ous function $K\,:\,\mathbb{R}_+\to\mathbb{R}_+$, such that
\begin{equation}\label{glestfin}
||\mathcal{N}({\bf v}(t))||_{\mathcal{X}}\le K(t),\quad \text{for}\quad t\in [0,t_1).
\end{equation}
Let $\mathcal{T}$ be the semigroup generated by $\mathcal{L}$. From the variation of constants formula, using Theorem 1.4.3 in \cite{Henry}, we have for $t\in[0,t_1)$ and for $\alpha\in \left(\frac{1}{2},1\right)$
\begin{align}
||{\bf v}(t)||_{\mathcal{X}^\alpha}\le & \left|\left|\mathcal{T}(t){\bf v}_0\right|\right|_{\mathcal{X}^\alpha}+\int_0^t\left|\left|\mathcal{T}(t-s)\mathcal{N}({\bf v}(s)) \right|\right|_{\mathcal{X}^\alpha}\,\ud s \nonumber \\
\le & \left|\left|\mathcal{T}(t){\bf v}_0\right|\right|_{\mathcal{X}^\alpha}+\int_0^t \left|\left|\mathcal{A}_1^\alpha\mathcal{T}(t-s)\right|\right|_{\mathcal{X}}||\mathcal{N}({\bf v}(s))||_{\mathcal{X}}\,\ud s \nonumber \\
\le & ||{\bf v}_0||_{\mathcal{X}^\alpha}Me^{at}+\int_0^t C_\alpha (t-s)^{-\alpha}e^{a(t-s)}K(s)\,\ud s.
\end{align}
Since $\alpha<1$, the integral is convergent, and we obtain the a priori bound
\begin{equation}
||{\bf v}(t)||_{\mathcal{X}^\alpha}\le H(t), \quad t\in [0,t_1),
\end{equation}
where $H\,:\,\mathbb{R}_+\to\mathbb{R}_+$ is continuous. This concludes the proof of the theorem. \eofproof

\begin{remark} 
Note that, to obtain the a priori upper bound in the proof of Theorem \ref{global}, we had to assume that $\Gamma\le 1$. 
On the other hand we proved local existence for any finite $\Gamma$. 
For $\Gamma>1$ finite time blow-up might occur, with $\int_0^1v(y,t)\,\ud y$ remaining bounded. 
To prove (or disprove) the finite time blow-up for $\Gamma>1$ is left as an open problem. 
We note that blow-up results were discussed for nonlinear reaction diffusion equations, where the $L^1$ norm of the solution is preserved, for example in \cite{Budd}. On the other hand, the blow-up phenomenon in \cite{Budd} is essentially possible due to the lack of positivity of solutions.

Here we only point out that a similar phenomena can be observed in the following ordinary differential equation
\begin{equation*}
\frac{\ud}{\ud t}v(t)=-\varrho v(t)+\beta v^p(t).
\end{equation*}
It is easily seen that for $p>1$ finite time blow-up occurs. Note that, of course, in the equation above the total population is not preserved. In contrast, in our model \eqref{superspread} the total population size remains constant. But we claim that the behaviour of the equation above might be mimicked in the structured model. For some fixed value of the structuring variable $y$, $v(y,t)$ could tend to infinity in finite time. Hence we may argue that the critical exponent could be $2$. This might be because in the structured model the diffusion has a distributing/smoothing effect.
\end{remark}

\section{Existence of (positive) steady states}

In this section we investigate the existence of steady states of model \eqref{superspread}. First we note that, a continuum family of ``semi-trivial'' steady states of the form $(0,S_*)$ exists (i.e. for any $S_*\ge 0$).
The interesting question is when, and under what assumptions may a positive (endemic) steady state of the model exist.
We will address this question by reformulating the steady state problem as a spectral problem for an appropriately defined
linear operator. First we discuss the classic bilinear model, that is when $\gamma\equiv 0$. In this case
we are able to prove in a relatively straightforward fashion the existence of endemic steady states.
Then we will consider the much more delicate case when we have a highly nonlinear infection process.

\subsection{Existence of the endemic steady state when $\gamma\equiv 0$}
In case of $\gamma\equiv 0$, our model \eqref{superspread} resembles the classic bilinear ODE model, with mass action assumption for the recruitment of newly infected individuals. Our main result below is in direct comparison with the corresponding result for the ODE model.
\begin{theorem}\label{steadyth0}
Assume that $\gamma\equiv 0$, $\varrho\not\equiv 0$, and $\int_0^1\beta(x,y)\,\ud y>0$, for every $x\in [0,1]$.
Then there exists a unique value $S_*$ of the susceptible population size, such that for any infected population size $V_*>0$,  there is a unique steady state $(v_*,S_*)$, with $V_*=\int_0^1 v_*(x)\,\ud x$.
\end{theorem}
\noindent {\bf Proof.}
For every $R\in\mathbb{R}_+$ we define the following linear operator
\begin{align}
\Psi_R\,v= & \left(d(\cdot)v'\right)'-\varrho(\cdot)v+R\int_0^1\beta(\cdot,y)v(y)\,\ud y, \label{operator1} \\
\text{D}\left(\Psi\right)= & \left\{v\in W^{2,1}(0,1)\,|\, v'(0)=v'(1)=0\right\}.\label{operator2}
\end{align}
For every $R\in\mathbb{R}_+$ the operator $\Psi_R$ generates a positive irreducible and analytic semigroup of linear operators on $L^1(0,1)$. Also note that the resolvent $R(\lambda,\Psi_R)$ is compact. Hence the spectral bound $s(\Psi_R)$ is a simple
eigenvalue with a corresponding strictly positive eigenvector. In fact, the spectral bound is the only eigenvalue with a corresponding strictly positive eigenvector. Therefore, $(v_*,S_*)\in\mathcal{X}^\alpha$ is a positive steady state of model
\eqref{superspread} (with $\gamma\equiv 0$) if and only if $s(\Psi_R)=0$ for some $R>0$,  $v_*\in\text{Ker}(\Psi_R)$, and $S_*=R$. Note that in this case, the second steady state equation,
\begin{equation}
\int_0^1\varrho(x)v_*(x)\,\ud x=R\int_0^1\int_0^1\beta(x,y)v_*(y)\,\ud y\,\ud x,\label{equation1}
\end{equation}
is automatically satisfied, as it is simply the integral of the equation $\Psi_R\, v=0$.

Next we show that  $s(\Psi_0)<0$. To this end we multiply the eigenvalue equation
\begin{equation*}
\Psi_0\,v=s\left(\Psi_0\right)\,v
\end{equation*}
by  $v$, and integrate from $0$ to $1$. We obtain
\begin{equation}\label{eigenvector}
-\int_0^1\left(d(x)(v'(x))^2+\varrho(x)v^2(x)\right)\,\ud x=s\left(\Psi_0\right)\int_0^1v^2(x)\,\ud x.
\end{equation}
Note that any eigenvector $v$ belongs to the domain of the operator $\Psi_0$, and we have $v\in W^{2,1}(0,1)\subset W^{1,2}(0,1)$, by the Sobolev embedding theorem. Hence equation \eqref{eigenvector} shows that the spectral bound $s\left(\Psi_0\right)$ is negative, if $\varrho\ge 0$ and $\varrho\not\equiv 0$.

To show that the function $f\,:\,R\to s\left(\Psi_R\right)$ is strictly monotone increasing, we use Proposition A.2 from \cite{AB}, which we stated after the proof of this theorem for the reader's convenience. In particular, in our setting we apply Proposition \ref{props} below with $\mathcal{A}_1=\Psi_{R_1}$ and $\mathcal{A}_2=\Psi_{R_2}$, where $R_1<R_2$, and $R_1,R_2\in\mathbb{R}_+$.

Next we show that the spectral bound $s\left(\Psi_R\right)$ changes sign indeed as $R$ increases.
To this end we note that if $\mathcal{A}$ is a generator of a strongly continuous positive semigroup on $L^1$, then we have
\begin{equation}\label{AGGcorollary}
s(\mathcal{A})\ge \sup_{\mu}\,\{\mu\in\mathbb{R}\,|\,\mathcal{A}f\ge \mu f, \,\, \text{for some}\,\, 0<f\in D(\mathcal{A})\}.
\end{equation}
Note that this result is stated in Corollary 1.14 B-II in \cite{AGG}, for generators of semigroups on $C(K)$, where $K$ is a compact topological space. It is easily seen (directly from the proof of Corollary 1.14 B-II in \cite{AGG}) that \eqref{AGGcorollary} holds for generators of semigroups on $L^1$, too. For any (fixed) $R\in\mathbb{R}_+$ we let $\mathcal{A}=\Psi_{R}$. Hence applying \eqref{AGGcorollary} it is enough to show that there exists an $0<f\in D(\Psi_R)$ such that
\begin{equation}\label{AGGestimate}
\left(d(\cdot)f'\right)'-(\varrho(\cdot)+\varepsilon)f+R\int_0^1\beta(\cdot,y)f(y)\,\ud y\ge 0
\end{equation}
holds for some $\varepsilon>0$, and for a large enough $R\in\mathbb{R}_+$. Note that for a positive constant function $f$
 inequality \eqref{AGGestimate} holds, if
\begin{equation}\label{AGGestimate2}
R\int_0^1\beta(\cdot,y)\,\ud y\ge r+\varepsilon,
\end{equation}
since $\varrho\le r$ (see \eqref{hypo}). Since
$\beta$ is a continuous function, \eqref{AGGestimate2} holds if
\begin{equation*}
R\ge \frac{r+\varepsilon}{\displaystyle\min_{x\in [0,1]}\int_0^1\beta(x,y)\,\ud y}.
\end{equation*}
Hence there exists a unique positive $R_* (=S_*)$, such that $s(\Psi_{R_*})=0$, with a family of corresponding strictly positive eigenvectors: $\kappa v_*,\,\kappa\in\mathbb{R}_+$. That is, we have shown that there exists a family of endemic steady states of the form $(\kappa v_*,S_*)$, $\kappa\in\mathbb{R}_+$.
\eofproof

\begin{proposition} (Proposition A.2, \cite{AB})\label{props} Let $\mathcal{A}_1,\mathcal{A}_2$ be resolvent positive operators, with dense domain such that
\begin{equation*}
0\ll R(\lambda,\mathcal{A}_1)\le R(\lambda,\mathcal{A}_2),\quad\text{for}\quad \lambda>\max\{s(\mathcal{A}_1),s(\mathcal{A}_2)\}.
\end{equation*}
If $\mathcal{A}_1\ne \mathcal{A}_2$ and $s(\mathcal{A}_1),s(\mathcal{A}_2)$ are poles of the resolvent of $\mathcal{A}_1,\mathcal{A}_2$, respectively, then $s(\mathcal{A}_1)<s(\mathcal{A}_2)$.
\end{proposition}

\begin{remark}
Note that our result is in direct comparison with the corresponding result for the ODE model \eqref{ODEmodel}, when $\gamma=0$. In particular, it is easily seen that the corresponding ODE model admits endemic steady states of the form $\left(I_*,\frac{\varrho}{\beta}\right)$, for any $I_*>0$.
\end{remark}

\subsection{Existence of the endemic steady state for $\gamma\equiv 1$}

Next we discuss the existence of a positive steady state for the more difficult and interesting case, when $\gamma\equiv 1$, that is when we have a highly nonlinear infection process. We note that the general framework we developed very
recently in \cite{CF2}, to treat the steady state problem for a large class of equations, does not apply to our model \eqref{superspread}, because it incorporates an infinite dimensional nonlinearity. But our nonlinearity is monotone, hence the approach we employ here is somewhat similar (but with a lot of additional technical difficulties, as we will see later) to the one we developed in \cite{CF}.

We emphasize that, as we will see later, the main difficulties to prove the existence of the endemic steady state
in case of a highly nonlinear infection process are the following. First, we need to recast the steady state problem as a spectral problem for well-defined linear operators. This is non-trivial, since as we have seen earlier, the nonlinearity cannot be defined on the natural state space $L^1$. Therefore we need to use a slightly unusual parameter space. At the same time, we need to uniformly control the spectral bound of the family of operators over the parameter set.
Secondly, to apply Schauder's fixed point theorem on a particular (convex) subset of the parameter space,
we need to prove that this subset of the parameter space will be mapped into a set which is contained in a compact set.
This is far from trivial, even impossible for some $\gamma$, because of the local nonlinearity.
To overcome all of the difficulties at the same time is very challenging. To illustrate the boundaries of the approach we employ, in our main result below we provide proofs of partial results for the most general choice of $\gamma$.
\begin{theorem}\label{steadyth}
Assume that $\gamma\equiv 1$ holds, $\varrho\not\equiv 0$, and that $\beta$ is strictly positive. Then, for every $S_*>0$, the structured model \eqref{superspread} admits a (strictly) positive  (endemic) steady state of the form  $(v_*,S_*)$.
\end{theorem}
\noindent {\bf Proof.} We introduce a somewhat unusual parameter set, namely we consider the set
\begin{equation}
C=\{0\le u\in W^{1,1}(0,1)\,|\, 0<||u||_{W^{1,1}(0,1)}\le 1 \}.
\end{equation}
For every $(u,R)\in C\times\mathbb{R}_+$ we define the following linear operators:
\begin{align}
\Psi^1\,v= & \left(d(\cdot)v'\right)'-\varrho(\cdot)v,\quad \text{D}\left(\Psi^1\right)=\left\{v\in W^{2,1}(0,1)\,|\, v'(0)=v'(1)=0\right\},\label{oper1} \\
\Psi^2_{(u,R)}\,v= & R\int_0^1\beta(\cdot,y)v(y)u(y)^{\gamma(y)}\,\ud y, \label{oper2} \\
\Psi_{(u,R)}= & \Psi^1+\Psi^2_{(u,R)},\quad \text{D}\left(\Psi_{(u,R)}\right)=\text{D}\left(\Psi^1\right).\label{oper}
\end{align}
Similarly as in Theorem \ref{steadyth0}, it is shown that the operators $\Psi_{(u,R)}$ generate positive, irreducible and analytic semigroups. Therefore $(v_*,S_*)\in\mathcal{X}^\alpha$ is a positive steady state of model
\eqref{superspread} if and only if $s(\Psi_{(v_*,S_*)})=0$, $v_*\in\text{Ker}(\Psi_{(v_*,S_*)})$. Note that in this case, the
second steady state equation, which reads
\begin{equation}
\int_0^1\varrho(x)v_*(x)\,\ud x=S_*\int_0^1\int_0^1\beta(x,y)(v_*(y))^{1+\gamma(y)}\,\ud y\,\ud x,\label{sseq2}
\end{equation}
will be automatically satisfied.

Therefore, the proof of the theorem consist of two parts. First, we
need to show that the kernel of $\Psi_{(u,R)}$ contains strictly positive vectors. Secondly, we need to
show  the existence of a fixed point of a nonlinear map defined on the set $C$.

The operator $\Psi^1$ generates a positive irreducible and analytic semigroup of linear operators on $L^1(0,1)$. Since for every $R\in\mathbb{R}_+$ and $u\in C$ the operator $\Psi^2_{(u,R)}$ is positive and bounded on $L^1(0,1)$, we have that $\Psi_{(u,R)}$ generates a positive analytic semigroup on $L^1(0,1)$. The semigroup is also immediately compact due to the boundedness of the interval $(0,1)$, and the Sobolev embedding theorem.

Since the semigroup $\mathcal{S}_1$ generated by $\Psi^1$ is positive and $\Psi^2_{(u,R)}$ is positive, we have that $0\le\mathcal{S}_1\le \mathcal{S}$, where $\mathcal{S}$ is the semigroup generated by $\Psi_{(u,R)}$.
Since any ideal $I$ of $L^1(0,1)$, which is invariant for $\mathcal{S}$ is also invariant for $\mathcal{S}_1$, the irreducibility of $\mathcal{S}_1$ implies the irreducibility of $\mathcal{S}$. Therefore we have that the spectral bound $s(\Psi_{(u,R)})$ is a simple eigenvalue.

For every $R\in \mathbb{R}_+$  the function $f_R\,:\, u\to
s(\Psi_{(u,R)})$ is a continuous function from $W^{1,1}_+$ to
$\mathbb{R}$, see \cite[IV-3.5]{K}. Note that we have already shown
in the proof of Theorem \ref{steadyth0} that for every
$R\in\mathbb{R}_+$ we have  $s(\Psi^1)=s(\Psi_{({\bf 0},R)})<0$. 
This, together with the continuous dependence of the
spectral bound on the parameter $u$, implies that there exists an
$r_*>0$ such that for $||u||_{W^{1,1}}\le r_*$ we still have
$s(\Psi_{({\bf u},R)})<0$. This $r_*$ of course depends on $R$, in
general.

For any $0\not\equiv\gamma\ge 0$ and $R\in\mathbb{R}_+$ the spectral
bound $s(\Psi_{(u,R)})$ is strictly monotone along positive rays of
the parameter space, i.e. for every $R\in\mathbb{R}_+$ the function
$g_R\,:\, \theta \in\mathbb{R}_+ \to
s\left(\Psi_{(\theta u,R)}\right)$ is strictly monotone increasing,
for any $u\in C$. This is established, similarly as in the proof of
Theorem \ref{steadyth0}, using Proposition \ref{props}. In
particular we can apply Proposition \ref{props} with
$\mathcal{A}_1=\Psi_{(\theta_1 u,R)}$ and
$\mathcal{A}_2=\Psi_{(\theta_2 u,R)}$, where $\theta_1<\theta_2$.

Next we show that, for every $R\in\mathbb{R}_+$,
the spectral bound $s\left(\Psi_{(u,R)}\right)$ changes sign along
positive rays intersecting $C$, whenever $0\le \gamma$, and $\gamma$ only vanishes
on a set of measure zero. To this end, fix $u\in C.$ Then, arguing exactly as in the proof of Theorem
\ref{steadyth0}, it is enough to show that, for $\theta$ sufficiently large,
\begin{equation}\label{AGGest2}
R\int_0^1\beta(\cdot,y)\theta^{\gamma(y)}u(y)^{\gamma(y)}\,\ud y\ge r+\varepsilon.
\end{equation}
Let $U\subseteq [0,1]$ be the essential support of the function $u$. Then, the (Lebesgue) measure $m(U)$ of $U$ is
positive. On the other hand, there exists a $\delta>0$ such that the set $V\subseteq [0,1]$ on which
$\gamma(y)\ge \delta$ has measure larger than $1-m(U)$. Then the closed set $W=U\cap V$ has positive measure, and for $\theta\ge 1$ we have
\begin{align}
R\int_0^1\beta(\cdot, y)\theta^{\gamma(y)}u^{\gamma(y)}(y)\,\ud y & \ge R\,\theta^\delta\int_W\beta(\cdot,y)u^{\gamma(y)}(y)\,\ud y \nonumber \\
& \ge R\,\theta^\delta \min_{x\in [0,1]}\int_W\beta(x,y)u^{\gamma(y)}(y)\,\ud y \to\infty, \label{spectrestimate}
\end{align}
as $\theta\to\infty$, since
\begin{equation}
\min_{x\in [0,1]}\int_W\beta(x,y)u^{\gamma(y)}(y)\,\ud y>0.
\end{equation}

So we have shown that for every $R\in\mathbb{R}_+$ the spectral bound function changes sign along rays intersecting
$C$. We now introduce the level set where the spectral bound vanishes. That is, for any $R\in\mathbb{R}_+$ we
define the set
\begin{equation}
S_R=\{0<u\in W_+^{1,1}(0,1)\,|\, s(\Psi_{(u,R)})=0\}.
\end{equation}
Note that, although the argument above proves that for any
$R\in\mathbb{R}_+$ the spectral bound changes sign along every ray
intersecting $C$, it does not suffice as we will need to control the norm of the level set $S_R$, too. To obtain the
boundedness of the level set $S_R$ even in $L^1$, we need a further
condition on $\gamma$. In particular for $\gamma\ge 1$ we can
control the growth behaviour of the spectral bound function
uniformly, and therefore we obtain an a-priori bound of the level
set $S_R$, as follows. Let us consider the following set
\begin{equation*}
C'=\left\{u\in W_+^{1,1}(0,1)\,|\,
||u||_1=1\right\}.
\end{equation*}
We are going to prove that, for any fixed $R$, every element of $S_R$ can be obtained
by multiplying an element of $C'$ with a constant which is uniformly bounded. To this end, we prove that for $u$ from
$C'$ the spectral bound of $\Psi_{(u,R)}$ is bounded below uniformly in $u$.

Let $u\in C'$ and let us denote by $\Omega$ the subset of $[0,1]$ on
which $u<1$. Assuming $1\le\gamma(y)\le\Gamma$ and
$\theta \ge 1$, we obtain the following estimate.
\begin{align}
\int_0^1\beta(\cdot,y)\theta^{\gamma(y)}u^{\gamma(y)}\,\ud y \ge &  \min_{x,y\in [0,1]}\{\beta(x,y)\}\,\theta\,\left(\int_{\Omega}u(y)^{\gamma(y)}\,\ud y+\int_{[0,1]\setminus\Omega}u(y)^{\gamma(y)}\,\ud y\right) \nonumber \\
\ge &  \min_{x,y\in [0,1]}\{\beta(x,y)\}\,\theta\,\left(\int_{\Omega}u(y)^\Gamma\,\ud y+\int_{[0,1]\setminus\Omega}u(y)\,\ud y\right) \nonumber \\
\ge & \min_{x,y\in[0,1]}\{\beta(x,y)\}\,\theta\,\left(\left(\int_{\Omega}u(y)\,\ud y\right)^\Gamma+\int_{[0,1]\setminus\Omega}u(y)\,\ud y\right) \nonumber \\
= & \min_{x,y\in [0,1]}\{\beta(x,y)\}\,\theta\,\left(\int_{[0,1]\setminus\Omega}u(y)\,\ud y+\left(1-\int_{[0,1]\setminus\Omega}u(y)\,\ud y\right)^\Gamma\right). \label{parestimate1}
\end{align}
The last inequality holds, since by H\"{o}lder's inequality
\begin{equation*}
\int_\Omega u(y)\,\ud y\le m(\Omega)\left(\int_\Omega
u(y)^\Gamma\,\ud y\right)^{\frac{1}{\Gamma}},
\end{equation*}
which implies that
\begin{equation*}
\int_\Omega u(y)^\Gamma\,\ud y\ge \left(\int_\Omega
u(y)\,\ud y\right)^\Gamma.
\end{equation*}
If $\Gamma=1$ then the last line in \eqref{parestimate1} simply reads $\displaystyle\min_{x,y\in [0,1]}\{\beta(x,y)\}\,\theta$.
For $\Gamma>1$, note that the function
\begin{equation*}
f_\Gamma(x)=x+(1-x)^\Gamma,\quad x\in [0,1],
\end{equation*}
has a unique minimum value of
\begin{equation*}
\Delta:=1-\Gamma^{\frac{1}{1-\Gamma}}+\Gamma^{\frac{\Gamma}{1-\Gamma}}>0,
\end{equation*}
on the interval $[0,1]$, attained at the point $1-\Gamma^{\frac{1}{1-\Gamma}}$. Hence there exists a $\Delta>0$, such that
\begin{equation}
\int_0^1\beta(\cdot,y)\theta^{\gamma(y)}u(y)^{\gamma(y)}\,\ud y\ge \displaystyle\min_{x,y\in [0,1]}\{\beta(x,y)\}\,\theta\,\Delta,
\end{equation}
uniformly for $u\in C'$. In particular, for any $R\in\mathbb{R}_+$ and $u\in C'$, we have that \eqref{AGGest2} holds, if
\begin{equation}
\theta\ge\frac{r+\varepsilon}{R\,\Delta\,\displaystyle\min_{x,y\in [0,1]}\{\beta(x,y)\}}.\label{thetaeq}
\end{equation}
This means that for every $R\in\mathbb{R}_+$ there exists a constant $\theta_*\in\mathbb{R}_+$, such that
$s(\Psi_{(\theta_*u,R)})>0$, for every $u\in C'$. That is, for every $s\in S_R$ we have that $||s||_1\le \theta_*$, that
is, the level set $S_R$ is bounded in $L^1$.

Note that the set $C$ is convex. Our goal is to apply Schauder's fixed point theorem, for every $R\in\mathbb{R}_+$, 
to an appropriately defined map $\Phi_R$, on the set $C$. In particular, we define a family of nonlinear maps
$\Phi_R\,:C\to C$, where $R\in\mathbb{R}_+$ such that the maps $\Phi_R$ are continuous and compact. Each $c\in C$ determines a unique positive ray $K_c:=\{\theta c \,|\, \theta\in\mathbb{R}_+\}$.
On any $C$-intersecting ray $K_c$, there is a unique element $u_*\in
S_R$, that is for which $s\left(\Psi_{(u_*,R)}\right)=0$, with a
corresponding strictly positive unique normalised (in
$W^{1,1}(0,1)$!) eigenvector $v_{u_*}\in C$. That is,  we have
\begin{equation}\label{steadyeq1}
\Psi_{(u_*,R)}v_{(u_*)}=\left(dv'_{(u_*)}\right)'-\varrho(\cdot)v_{(u_*)}+R\,\int_0^1\beta(\cdot,y)v_{(u_*)}(y)u_*(y)^{\gamma(y)}\,\ud y = 0.
\end{equation}
Hence for every $R\in\mathbb{R}_+$ we define
\begin{equation}
\Phi_R\,:\underbrace{c}_{\in C}\to\, \underbrace{u_*}_{\in K_c\cap S_R}\,\to\,\underbrace{v_{u_*}}_{W^{2,1}(0,1)\cap C},
\end{equation}
that is, $\Phi_R(c)=v_{u_*}$.
The maps $\Phi_R$ are continuous, because the projection along rays is continuous on $C$, and the function $u_*\to v_{u_*}$ is also continuous. This is contained in \cite{K}.

For compactness, note that $W^{2,1}(0,1)$ is compactly embedded in
$W^{1,1}(0,1)$. It is left to prove that the set of eigenvectors 
$\Phi_R(C)$ is bounded in $W^{2,1}(0,1)$. Let
$v=v_{u_*}$ be the strictly positive normalised eigenvector of
$\Psi_{(u_*,R)}$ corresponding to its spectral
bound. For any $u_*\in S_R$, and assuming $\gamma\le 1$, we have
\begin{align}
||v||_{W^{2,1}} & =||v||_1+||v'||_1+||v''||_1\le 1 +\frac{1}{d_0}\left|\left|\Psi^1\,v-d'v'+\varrho v\right|\right|_1 \nonumber \\
& \le \frac{1}{d_0}\left|\left|\Psi^1\,v\right|\right|_1+\max\left(\frac{d_1}{d_0},\frac{r}{d_0}\right)||v||_{W^{1,1}} + 1 \nonumber \\
& =\frac{1}{d_0}\left|\left|\Psi^2_{(u_*,R)}\, v\right|\right|_1+D \nonumber \\
& \le R\,b\,||v||_\infty\int_0^1 u_*^{\gamma(y)}(y)\,\ud y +D\le R\, \hat{b}\,\left(1+\int_0^1 u_*(y)\,\ud y \right)+D \nonumber \\
& \le R\,\hat{b}(1+\theta_*)+D, \label{finalestimates}
\end{align}
for some constants $\hat{b},D$. Therefore, $\overline{\Phi_R(C)}$ is
a compact set contained in the unit sphere of $W^{1,1}(0,1)$, and in the positive cone, and so in $C$. Hence
Schauder's fixed point theorem implies the existence of a fixed point of $\Phi_R$, which we denote by $v_*^R$. This means that from
equation \eqref{steadyeq1} we have
\begin{equation}
0=\left(d\left(v_{*}^{R}\right)'\right)'-\varrho(\cdot)v_*^R+R\,\int_0^1\beta(\cdot,y)v_*^R(y)\left(\kappa\,v_*^R(y)\right)^{\gamma(y)}\,\ud y,\label{steadyeq2}
\end{equation}
for some $\kappa>0$. Hence the positive steady state of \eqref{superspread}  is
\begin{equation*}
\left(\kappa\,v_*^R,R\right)=:(v_*,S_*).
\end{equation*}
Note that since the eigenvector $v_*^R$ is normalised in $W^{1,1}$, we have that the total po\-pu\-lation size $\kappa\int_0^1v_*^R(x)\,\ud x+R$ of the steady state $(v_*,S_*)$ is bounded above by $\kappa+R$.
\hfill $\Box$

\begin{remark}
We note that since $\kappa v_*^R\in S_R$, we have in fact that the total population size of the infected individuals,
which equals $||\kappa v_*^R||_1$, is bounded by $\theta_*$. At the same time we note that $\theta_*$ depends on the growth rate at which the spectral bound function $g_R\,:\, \theta \,\to\, s\left(\Psi_{(\theta u,R)}\right)$ increases
along positive rays in the parameter space. In particular, we can see from the proof of Theorem \ref{steadyth}, and from inequality \eqref{thetaeq}, that $\theta_*$, the upper bound for the $L^1$ norm of the level set $S_R$, only has to satisfy
\begin{equation*}
\frac{\displaystyle\max_{x\in [0,1]} \varrho(x)}{R\,\Delta\displaystyle\min_{x,y\in [0,1]}\beta(x,y)}\le\theta_*.
\end{equation*}
This implies that there exists a continuous monotone decreasing function $\tau$ of $R$, such that
\begin{equation*}
\displaystyle\lim_{R\to\infty}\tau(R)=0,\quad \text{and}\quad ||S_R||_1\le\tau(R).
\end{equation*}
Hence the total population size of the endemic steady state is bounded above by $\tau(R)+R$.
\end{remark}

\begin{remark}
We note that the result above is comparable to the corresponding result for the unstructured model \eqref{ODEmodel}.
The unstructured model admits the steady state $(0,S_*)$ for any $S_*\ge 0$, and any positive steady state
$(V_*,S_*)$ with total population size $P_*=V_*+S_*$ satisfies the equation
$$\frac{\varrho}{\beta}=S_*V_*^\gamma.$$
It is shown (from $P=V+V^{-\gamma}\frac{\varrho}{\beta}$) that the critical population size at which positive equilibria emerge
is
\begin{equation*}
\bar{P}=\left(\gamma\frac{\varrho}{\beta}\right)^{\frac{1}{1+\gamma}}+\frac{\varrho}{\beta}\left(\gamma\frac{\varrho}{\beta}\right)^{-\frac{\gamma}{1+\gamma}}=\left(\gamma\frac{\varrho}{\beta}\right)^{\frac{1}{1+\gamma}}\left(1+\gamma^{-1}\right).
\end{equation*}
In particular for population sizes $P_*$ greater than $\bar{P}$
there are two positive equilibria. This seems to be different than
for the structured model \eqref{superspread}, for which we showed in
the previous theorem that there is (at least) one
positive equilibrium. Notice however that we proved the existence of
a positive equilibrium for any fixed value $R$, which is the
population size of the susceptible individuals. Indeed, it is the
case already in the unstructured model that for any given
susceptible population size $S$ there is a unique positive steady
state $\left(S,\left(\frac{\varrho}{\beta S}\right)^{\frac{1}{\gamma}}\right)$.
\end{remark}

\begin{remark}
As we can see from the proof of Theorem \ref{steadyth} we could only obtain an a-priori uniform $L^1$ bound for the level set $S_R$, when $\gamma\ge 1$, and we cannot obtain an a-priori $W^{1,1}$ bound for the level set $S_R$ for any values of $\gamma$.
At the same time, to use the $L^1$ bound of the level set $S_R$ to prove that the image $\Phi_R(C)$ is contained in a compact subset of $C$ was only possible for $\gamma\le 1$.
\end{remark}

\section{Concluding remarks}

In this paper we considered an SIS-type partial differential equation model for the spread of an infectious disease.
Our goal was to introduce a model, which takes into account a super-spreading phenomenon, that is when individuals
exhibit (typically higher than normal) levels of infectiousness. Modelling this interesting phenomena resulted in a highly nonlinear recruitment term (infection incidence). In the context of ordinary differential equations, this (superspreading)  phenomenon was already investigated in the seminal paper \cite{CS}, and more recently for example in \cite{KOR}. In a structured partial differential equation model the local nonlinearity leads to substantial difficulties. In particular we expect that blow-up phenomena may occur, even though the total population size is preserved. To prove (or disprove) this claim is left for future work. 

Since the model we developed incorporated a nonlinearity, which could not be defined on the biologically natural state space, we used the framework of fractional power spaces as developed in \cite{Henry}, to prove global existence of solutions up to a critical value of the exponent. This method utilizes the (advantageous) fact that the linear part of the problem
is governed by an analytic semigroup. An added advantage of the semilinear formulation developed in \cite{Henry} is that the Principle of Linearised Stability can be immediately established.

We also studied the existence of the endemic steady state. Since the steady state equation is a second order integro-differential equation we used an implicit approach to prove the existence of the endemic steady state. The main idea
is to reformulate the steady state problem as a family of abstract eigenvalue problems and to cast a fixed point problem.
A similar approach was used recently in \cite{CF}, see also \cite{CF2} for a general theory for models with two-dimensional nonlinearities. In contrast to the model we treated in \cite{CF}, the local, but infinite dimensional nonlinearity in model \eqref{superspread} proved to be (not surprisingly) extremely challenging. In particular we had to choose an unnatural parameter space to work with, that is $W^{1,1}$, in which the unit sphere is not convex. Furthermore, any bounded subset of the parameter space contains elements with arbitrarily small $L^1$ norm, which means that the growth rate of the spectral bound
of the operator $\Psi$, and in turn the $W^{1,1}$ norm of its level set $S_R$, could not be controlled, except in some particular cases. This then resulted in a serious restriction when proving the compactness of the fixed point map.
In particular, all of this meant that we could only establish the existence of the endemic steady state
for the bilinear model and for the case of a quadratic nonlinearity. This is somewhat unsatisfactory as one may naturally expect
that strictly positive steady states may exist for functions $\gamma$ taking values between $0$ and $1$, too. Partly for this reason we gave proofs of the partial results in Theorem \ref{steadyth} to illustrate the boundaries of the method for model
\eqref{superspread}. All of these challenges underline the difficulty of capturing the super-spreading phenomenon in an 
infinite dimensional dynamical system.

We note that the next natural step would be to consider the stability of the steady states of model \eqref{superspread}.
As we noted before, the semilinear formulation readily allows to establish the Principle of Linearised Stability. The stability results, which can be 
directly applied to our model, are formulated in Theorems 5.1.1 and 5.1.3 and in Corollary 5.1.6 in \cite{Henry}.

In particular, for any $0\le \gamma$, the linearisation of equation \eqref{superspread} around the steady state $(v_*,S_*)$ reads

\begin{align}
w_t(x,t)= & (d(x)w_x(x))_x-\varrho(x)w+\int_0^1\beta(x,y)v_*(y)^{\gamma(y)}\left[R(t)v_*(y)+S_*w(y)(1+\gamma(y))\right]\,\ud y, \nonumber \\
R'(t)= & \int_0^1\varrho(x)w(x)\,\ud x-\int_0^1\int_0^1\beta(x,y)v_*(y)^{\gamma(y)}\left[R(t)v_*(y)+S_*w(y)(1+\gamma(y))\right]\,\ud y\,\ud x,\label{linearised}
\end{align}
with the appropriate boundary conditions. The linearised system \eqref{linearised} is governed by a compact analytic semigroup (since $v_*\in W^{2,1}(0,1)$), hence stability is governed by the eigenvalues of the generator of the linearised semigroup. The linearised equations \eqref{linearised} lead to the following eigenvalue problem

\begin{align}
\lambda w=\left(d(\cdot)w'\right)'-\varrho(\cdot)w &+\int_0^1\beta(\cdot,y)\left[Rv_*(y)^{1+\gamma(y)}+S_*v_*(y)^{\gamma(y)}(1+\gamma(y))v(y)\right]\,\ud y,\label{evalue1} \\
\lambda R=\int_0^1\varrho(x)w(x)\,\ud x &-R\int_0^1\int_0^1\beta(x,y)v_*(y)^{1+\gamma(y)}\,\ud y\,\ud x \nonumber \\
& -S_*\int_0^1\int_0^1\beta(x,y)v_*(y)^{\gamma(y)}(1+\gamma(y))w(y)\,\ud y\,\ud x. \label{evalue2}
\end{align}
It is outside the scope of the current paper to try to present a complete study of the eigenvalue problem above. However, 
the stability of the semi-trivial (disease free) steady state can be discussed relatively easily, and it may give an insight into the dynamic behaviour of the whole system.

For $\gamma\equiv 0$ the eigenvalue problem \eqref{evalue1}-\eqref{evalue2} at the steady state $(0,S_*)$ simply reads
\begin{equation}\label{evalueeq1}
\lambda\,w=\Psi_{S_*}\, w,\quad \lambda\, R=-\int_0^1\left(\Psi_{S_*}\,w\right)(x)\,\ud x,
\end{equation}
where $\Psi$ is defined in \eqref{operator1}. Integrating the first equation and adding it to the second one yields
\begin{equation}
\lambda\left(R+\int_0^1 w(x)\,\ud x\right)=0,
\end{equation}
which shows that $\lambda=0$ or $R=-\int_0^1 w(x)\,\ud x$, which implies that the second equation in \eqref{evalueeq1} is simply the integral of the first one. Note that the eigenvalue $\lambda=0$ (with eigenvector $(0,1)$) arises due to the fact that there is a continuum of equilibria, that is for any positive population size $S_*$. We have proven in Theorem \ref{steadyth0} that in fact there exists a unique critical value $S_*$ such that the spectral bound of $\Psi$ changes sign, i.e. it becomes positive, that is the steady state $(0,S_*)$ becomes unstable. This behaviour is in accordance with the corresponding result for the ordinary differential equation model \eqref{ODEmodel}.

It is also easily seen that for  $\gamma\not\equiv 0$ and $\varrho\not\equiv 0$, the steady state $(0,S_*)$ is locally asymptotically stable, within the level set $||v||_1+|S|=S_*$.
For $\gamma\not\equiv 0$ at the steady state $(0,S_*)$, the eigenvalue problem \eqref{evalue1}-\eqref{evalue2} reduces to
\begin{equation}
\lambda v=\left(d(\cdot)v'\right)'-\varrho(\cdot)v,\quad v'(0)=v'(1)=0,\quad  \lambda S=\int_0^1\varrho(x)v(x)\,\ud x.\label{evspec}
\end{equation}
Multiplying the first equation in \eqref{evspec} by $v$ and integrating from $0$ to $1$ we obtain
\begin{align*}
\lambda\int_0^1 v^2(x)\,\ud x=-\int_0^1\varrho(x)v^2(x)\,\ud x,
\end{align*}
which shows that any eigenvalue must be negative if $\varrho\not\equiv 0,\, v\not\equiv 0$. If $v\equiv 0$, the last equation in
\eqref{evspec} implies that $\lambda=0$ is an eigenvalue with (normalised) eigenvector $(0,1)$, which corresponds to the fact that $(0,S_*)$, for  $S_*\neq 0$, is a continuum of disease free equilibria. This result again, is in accordance with the
corresponding result for the unstructured ordinary differential equation model \eqref{ODEmodel}.

Finally we note that it will be an interesting question to address in future work, 
how to incorporate more compartments, such as one for the recovered class (R), and also population dynamics (i.e. birth and death processes), in model \eqref{superspread}, at least for the case $\gamma\le 1$.

\section*{Acknowledgments}
\`{A}. Calsina was partially supported by the research projects 2009SGR-345 and DGI MTM2011-27739-C04-02,
and by the Edinburgh Mathematical Society while visiting the University of Stirling.
J. Z. Farkas was supported by a research grant from the The Carnegie Trust, and by the research project DGI MTM2011-27739-C04-02, while visiting the Universitat Aut\`{o}noma de Barcelona.

\end{document}